
\input gtmacros
\input amsnames
\input amstex
          
\catcode`\@=12        

\input gtmonout
\volumenumber{2}
\volumeyear{1999}
\volumename{Proceedings of the Kirbyfest}
\pagenumbers{157}{175}
\papernumber{9}
\received{23 February 1999}\revised{26 May 1999}
\published{18 November 1999}

\let\\\par

\let\gttitle\title
\def\title#1\endtitle{\gttitle{#1}}
\let\gtauthor\author
\def\author#1\endauthor{\gtauthor{#1}}
\let\gtaddress\address
\def\address#1\endaddress{\gtaddress{#1}}
\let\gtemail\email
\def\email#1\endemail{\gtemail{#1}}
\def\subjclass#1\endsubjclass{\primaryclass{#1}}
\let\gtkeywords\keywords
\def\keywords#1\endkeywords{\gtkeywords{#1}}
\def\heading#1\endheading{{\def\S##1{\relax}\def\\{\relax\ignorespaces}
    \section{#1}}}
\def\head#1\endhead{\heading#1\endheading}

\def\subhead#1\endsubhead{\sh{#1}}
\def\subsubhead#1\endsubsubhead{\sh{#1}}
\def\specialhead#1\endspecialhead{\sh{#1}}
\def\demo#1{\rk{#1}\ignorespaces}
\def\enddemo{\ppar}

\def\qed{\ifmmode\quad\sq\else\hbox{}\hfill$\sq$\par\goodbreak\rm\fi}  
\def\proclaim#1{\rk{#1}\sl\ignorespaces}
\def\endproclaim{\rm\ppar}
\def\cite#1{[#1]}
\newcount\itemnumber

\let\itemold\item
\def\item{\itemold{{\rm(\number\itemnumber)}}%
\global\advance\itemnumber by 1\ignorespaces}
\def\S{section~\ignorespaces}  
\def\date#1\enddate{\relax}
\def\thanks#1\endthanks{\relax}   
\def\dedicatory#1\enddedicatory{\relax}  
\let\footnote\plainfootnote

\def\Refs{\ppar{\large\bf References}\ppar\bgroup\leftskip=25pt
\frenchspacing\parskip=3pt plus2pt\small}       
\def\endRefs{\egroup}
\def\widestnumber#1#2{\relax}
\def\endrefitem{}
\def\refdef#1#2#3{\def#1{\leavevmode\unskip\endrefitem#2\def\endrefitem{#3}}}
\def\ref{\par}
\def\endref{\endrefitem\par\def\endrefitem{}}
\refdef\key{\noindent\llap\bgroup[}{]\ \ \egroup}
\refdef\no{\noindent\llap\bgroup[}{]\ \ \egroup}
\refdef\by{\bf}{\rm, }
\refdef\manyby{\bf}{\rm, }
\refdef\paper{\it}{\rm, }
\refdef\book{\it}{\rm, }
\refdef\jour{}{ }
\refdef\vol{}{ }
\refdef\yr{$(}{)$ }
\refdef\ed{(}{ Editor) }
\refdef\publ{}{ }
\refdef\moreref{}{ }
\refdef\inbook{from: ``}{'', }
\refdef\pages{}{ }
\refdef\page{}{ }
\refdef\paperinfo{}{ }
\refdef\bookinfo{}{ }
\refdef\publaddr{}{ }
\refdef\eds{(}{ Editors)}
\refdef\bysame{\hbox to 3 em{\hrulefill}\thinspace,}{ }
\refdef\toappear{(to appear)}{ }
\refdef\issue{no.\ }{ }

\newcount\refnumber\refnumber=1
\def\refkey#1{\expandafter\xdef\csname cite#1\endcsname{\number\refnumber}%
\global\advance\refnumber by 1}
\def\cite#1{[\csname cite#1\endcsname]}
\def\Cite#1{\csname cite#1\endcsname}  
\def\key#1{\noindent\llap{[\csname cite#1\endcsname]\ \ }}

\refkey {BE}
\refkey {BHMV}
\refkey {G1}
\refkey {G2}
\refkey {H}
\refkey {KL}
\refkey {KM1}
\refkey {KM2}
\refkey {L}
\refkey {MR}
\refkey {M}
\refkey {Mu}
\refkey {RT}
\refkey {R}
\refkey {T}
\refkey {TV}
\refkey {Wa}
\refkey {W}

\def\demon{\demo{\bf Proof}}

\def\a{\alpha}

\def\b{\beta}

\def\g{\gamma}

\def\d{\delta}

\def\D{\Delta}

\def\z{\zeta}

\def\k{\kappa}

\def\l{\lambda}

\def\x{\xi}

\def\si{\sigma}

\def\t{\tau}

\def\o{\omega}

\def\O{\Omega}

\def\bz{\Bbb Z}
\def\bc{\Bbb C}
\def\({\left(}
\def\){\right)}
\def\]{\right]}
\def\[{\left[}
\def\slash{/}

\title  Quantum
invariants of periodic three-manifolds \endtitle

\author Patrick M Gilmer \endauthor \address 
Louisiana State University\\Department of Mathematics\\
Baton Rouge, LA 70803, USA\endaddress
\email gilmer@math.lsu.edu\endemail

\abstract
Let $p$ be an odd prime and $r$ be
relatively prime to $p$. Let $G$ be a finite $p$--group.
Suppose an oriented 3--manifold $\tilde M$ has a free
$G$--action
with orbit space $M$.  We consider certain
Witten--Reshetikhin--Turaev $SU(2)$ invariants $w_r( M)$  in
$\bz[ \frac 1 {2r}, e^{\frac {2 \pi i}{8r}}]$. We will show that
$w_r(\tilde M)$ $\equiv  \k^{3\ \text{\rm def }(\tilde M \rightarrow
M)}(w_r(M))^{|G|} \pmod{p}$. Here
$\k=e^{\frac {2 \pi i (r-2)}{8r}}$, $\text{\rm def }$ denotes the
signature defect, and $|G|$ is the number of elements in $G$. We also give a
version of this result if $M$ and
$\tilde M$ contain framed links or colored fat graphs.  We give similar
formulas for non-free actions which hold for a specified finite set of values
for $r$.
\endabstract
\asciiabstract{Let p be an odd prime and r be relatively prime to p.  
Let G be a finite p-group.  Suppose an oriented 3-manifold M-tilde has
a free G-action with orbit space M.  We consider certain
Witten-Reshetikhin-Turaev SU(2) invariants w_r(M).  We will give a
fomula for w_r(M) in terms of the defect of M-tilde --> M and the
number of elements in G.  We also give a version of this result if M
and M-tilde contain framed links or colored fat graphs.  We give
similar formulas for non-free actions which hold for a specified
finite set of values for r.}

\primaryclass{57M10}\secondaryclass{57M12}

\keywords
$p$--group action,  lens space, quantum invariant,
Turaev--Viro invariant, branched cover, Jones polynomial, Arf invariant
\endkeywords
\asciikeywords{p-group action,  lens space, quantum invariant,
Turaev-Viro invariant, branched cover, Jones polynomial, Arf invariant}

\makeshorttitle

\centerline{\small \it Dedicated to Rob Kirby on his sixtieth birthday}

\section{Introduction}

Assume $p$ is an odd prime, and that $r$ is relatively prime to $p$,
and $r \ge 3$. Let $G$ be a finite $p$--group, with $|G|$ elements.
We let $\xi_{a}$ denote $e^{\frac {2 \pi i}{a}}$.

Let $ M$  be an oriented closed 3--manifold with an
embedded
$r$--admissibly colored fat trivalent
graph $J$.   We include the case that $J$ is empty. We consider the
Witten--Reshetikhin--Turaev
$SU(2)$ invariants
$w_r(M,J)\in R_{r}=\bz[ \frac 1 {2r}, \xi_{t}]$
[\Cite{W}, \Cite{RT}] where $t$ is $4r$ if $r$ is even and $8r$ if $r$ is
odd.  Here $w_r(M,J)$ is the version of the WRT--invariant which is
denoted $I_{-\xi_{4r}}(M,J)$ by Lickorish \cite{L}, assuming $\mu$ (also known
as $\eta$ in this paper ) is
chosen to be  $ \(\xi_{2r}-\xi_{2r}^{-1}\)/{(i\sqrt{2r})}$.
In terms of the Kirby--Melvin
\cite{KM1} normalization $\t_r$, one has
$w_r(M)=  \t_r(M)\ / \t_r(S^1 \times S^2)$.

Let $\tilde M$  be an oriented closed 3--manifold with a
$G$ action. The singular set of the  action is the collection of points
whose isotropy subgroup is non-trivial.
Let $\tilde J$ be an equivariant $r$--admissibly colored fat graph $\tilde J$ in
$\tilde M$ which is disjoint
from the singular set. The action is assumed to preserve the coloring and
thickening of $\tilde J$.  Burnside's theorem asserts that the center of a
finite $p$--group is non-trivial. It follows that the quotient map
$\tilde M \rightarrow \tilde M/G$ can be factored as a sequence of $\bz_p$
(possibly branched) covering maps. It follows that the orbit space is also a
closed 3--manifold with an $r$--admissibly colored fat graph, the image of
$\tilde J$. In this case, we will denote the orbit space of $\tilde M$
by $M$, and the orbit space of $\tilde J$ by $J$. We find
a relationship between $w_r(M,J)$, and $ w_r(\tilde M, \tilde J) $.

When the action is free, a signature defect of $\tilde M
\rightarrow M$ may be  defined as follows.  One can arrange that some number,
say $n$, of  disjoint
copies of $\tilde M \rightarrow M$ form the boundary of a regular $G$
covering space of a 4--manifold $W$, denoted by $\tilde W$. Then define
$$\text{\rm def }(\tilde M \rightarrow
  M) =  \frac 1 n \( {|G|}\ \text{Sign}(W)-    \ \text{Sign}(\tilde W)\).$$
The defect can be seen to be well-defined using
Novikov additivity together with that fact the signature of an $m$--fold
unbranched covering space of a closed manifold is $m$ times the signature of
the base manifold. This  generalizes the definition  of the signature defect
for a finite cyclic group
[\Cite {H}, \Cite{KM}]. As $\O_3(BG)=H_3(BG)$, and $H_3(BG)$
is annihilated by multiplication by $|G|$, $n$ can be taken to be $|G|$ in the
above definition. For this definition, it is not necessary that $G$ be a
$p$--group. We remark that $3\ \text{\rm def }(\tilde M \rightarrow M)$
is an integer. We will give a proof in \S 3.

We will be  working with congruences modulo the odd prime
$p$ in the ring of all algebraic integers (over $\bz$) after we have
inverted $2r$, where
$r$ is relatively prime to $p$.
Let $\k$ denote $\xi_{8} \xi_{4r}^{-1}$.

\proclaim{Theorem 1} If $G$ acts freely, then
$$ w_r(\tilde M,\tilde J)\equiv  \k^{3 \text{\rm def }(\tilde M
\rightarrow
M)}(w_r(M,J))^{|G|} \pmod{p}.$$ \endproclaim

Our equations in Theorems 1, 2 and 3 take place in $R_r \slash pR_r$.
We may think of $R_r$ as polynomials in $\x_t$ with coefficients in $\frac 1
{2r} \bz$ of degree less than $\phi(t)$. When multiplying such polynomial, one
should use the $t$-th  cyclotomic polynomial to rewrite the product as a
polynomial in degree less than $\phi(t)$.
Such a polynomial lies in $p R_r$ if and only if each of its coefficients maps
to zero
under the map $\frac 1 {2r} \bz \rightarrow \bz_p$ given by reduction modulo
$p$. For example we consider Theorem 1 with $r=5$ applied to the free $Z_3$
action on $S^3$ with quotient $L(3,1)$. This action has defect 2/3. In this case
$t=40$, and  the cyclotomic polynomial tells us to reduce to polynomials of
degree less than 16 via ${\x_t}^{16}= {\x_t}^{12}-{\x_t}^{8} +{\x_t}^{4}-1$. We
have:
$$w_5(L(3,1))= {\frac{3\,\x_t + 2\,{\x_t^3} - {\x_t^5} -
     4\,{\x_t^7} + 4\,{\x_t^9} +
     {\x_t^{11}} - 2\,{\x_t^{13}} -
     3\,{\x_t^{15}}}{10}},$$
$$w_5(S^3)= \eta ={\frac{\x_t + {\x_t^3} - 2\,{\x_t^5} +
     3\,{\x_t^7} + 3\,{\x_t^9} -
     2\,{\x_t^{11}} - 4\,{\x_t^{13}} +
     {\x_t^{15}}}{10}}$$
$$ w_5(S^3)-\k^2 (w_5(L(3,1)))^3 = 3\  {\frac{
      -\,{\x_t^5} + 2\,{\x_t^7} +
       2\,{\x_t^9} - 2\,{\x_t^{11}} -
       2\,{\x_t^{13}} + {\x_t^{15}}
        }{20}}$$
We note that $R_r \slash pR_r =
\bz[ \x_t] \slash p \bz[ \x_t]$.
Let $f$ be the smallest positive integer such that $p^f \equiv 1
\pmod{t}$.
$\bz[ \x_t] \slash p \bz[ \x_t]$ is the direct sum of $\phi(t) \slash
f$ fields each with
$p^f$ elements [\Cite{Wa}; page 14]. It is important that the $p$th power
 map is the Frobenius automorphism of $R_r \slash p R_r$.  So given $w_r(\tilde
M,\tilde J)$  and
$\text{\rm def }(\tilde M \rightarrow M)$, we can always solve uniquely  for
$w_r( M, J)$.
Theorem 1 by itself will provide no obstruction to the existence of a free
$G$--action on a given
manifold $\tilde M$.

If  $w_r(\tilde M, \tilde J) \ne 0$ for an infinite collection of $r$ prime to
$p$ then  the values of $w_r(\tilde M, \tilde J)$, and $w_r( M,  J)$ for this
collection of $r$ determine $\text{\rm def }(\tilde M \rightarrow M)$.
This is not apriori clear.

When the action is not free we have to restrict $r$  to a few values,
and we don't know an independent definition of the exponent of $\k$.
Theorems 2 and 3 by
themselves will provide no obstruction to  the existence of a
$G$--action on a given
manifold $\tilde M$.

\proclaim{Theorem 2} If $r$ divides $\frac{p \pm1}2$, then for some integer
$\d$,  one has
$$ w_r(\tilde M, \tilde J)\equiv  \k^{\d}(w_r(M,J))^{|G|} \pmod{p}.$$
If $G$ is cyclic and acts semifreely and $r$ divides $\frac{{|G|} \pm1}2$,
the same conclusion holds.  \endproclaim

We  tie down the factor $\k^{\d}$ in a special case.

\proclaim{Theorem 3} Suppose $r$ divides $\frac{p^s \pm1}2$, and $\tilde M$ is
a $p^s$--fold branched cyclic cover of a knot $K$ in a homology sphere
$M$. Let
$J$ be a colored fat graph in $M$ which misses $K$, and $\tilde J$ be
the inverse image of $J$ in $\tilde M$.
$$ w_r(\tilde M, \tilde J)\equiv \mp {\binom {-2r} p}^s
\k^{-3 \si_{p^s}(K)}(w_r(M,J))^{p^s}
\pmod{p},$$
\endproclaim

Here $\si_{p^s}(K)$ denotes the total ${p^s}$--signature of $K$ \cite{KM2}.
$\binom {-2r} p$ is a Legendre symbol.

As a corollary, we obtain the following generalization of a result of
Murasugi's [\Cite{M}; Proposition 8]. Murasugi's  hypothesis is that $\tilde L$ is
 $p$--periodic,  $\z$ is a { \it primitive}  $\frac {p \pm 1} 2$ th root of
unity and $p \ge 5$.
Here $V_L(\z)$ denotes the Jones polynomial evaluated at $t=\z$.

\proclaim {Corollary 1} Let $\tilde L$ be a $p^s$--periodic link in $S^3$
with quotient link $L$.
If $\z^{\frac {p^s \pm 1} 2}=1$ and $\z \ne -1$ , then
$$V_{\tilde L}(\z) \equiv V_{L}(\z^{\mp 1}) \pmod{p}.$$
\endproclaim

$L$ has an even number of components if and only if $\tilde L$ does also.
In this case, we must choose the same $\sqrt{\z}$ when evaluating both sides of
the above equation.
The above equation is  false, in general, for $\z=-1$:  the trefoil
is has period three  with orbit knot the unknot.
Making use of H. Murakami's formula \cite{Mu} relating the Jones polynomial
evaluated at $i$
to the Arf invariant of a proper link we have the following corollary. First we
observe: if $\tilde L$ is a $n$--periodic link in $S^3$
with quotient link $L$, and $n$ is odd, then $\tilde L$ is proper if and only
if
$L$ is proper.

\proclaim {Corollary 2} Let $\tilde L$ be a $n$--periodic proper link in $S^3$
with quotient link $L$. Let $n=\prod p_i^{s_i}$ be the prime factorization of
$n$.  Suppose for each $i$, $p_i^{s_i}\equiv \pm 1 \pmod{8}$ ( for each $i$ one
may choose $\pm $ differently), then $$\text{\rm Arf}(L) \equiv
\text{\rm Arf}(\tilde L) \pmod{2}.$$
\endproclaim

The period three action on the trefoil also shows the necessity of the
condition that $p_i^{s_i}\equiv \pm 1 \pmod{8}$.

In \S 2 we  establish  versions of Theorems 1 and 2 for the related
Turaev--Viro invariants by adapting an argument which Murasugi used to study the
bracket polynomial of periodic links \cite{M}. See also Traczyk's paper
\cite{T}.
In fact  these theorems ( Theorems 4 and 5) are immediate corollaries
of Theorems 1 and 2
but we prefer to give them direct proofs. The reason is that these
proofs
are simpler than the proofs of Theorems 1 and 2. These proofs can be
used to
obtain analogous results for other invariants
defined by Turaev--Viro type state sums.
Also  Lemma 3, that we establish to prove Theorem 5, is used later in
the
proofs of Theorems 2 and 3.

In \S 3, we relate $w_r(M,J)$ to the TQFT defined in \cite{BHMV}.
We discuss $p_{1}$--structures.  We also rephrase
Theorem 1 in terms of manifolds with $p_{1}$-- structure, and reduce the proof of
Theorem 1 to the case $G=\bz_p$.
We say that a regular $\bz_{p^s}$--cover which is a quotient of a regular
$\bz$ cover is a simple  $\bz_{p^s}$--cover.
In \S 4, we derive Theorem 1 for simple $\bz_{p^s}$--covers of closed manifolds.
This part of the argument applies  generally to  quantum invariants associated
to any TQFT.
We also obtain a version of Theorem 1 for simple $\bz_{p^s}$--covers of
manifolds whose
boundary is a torus.
In \S 5, we derive Theorem 1 in the case $M$ is a lens space and $G= \bz_p$.
In \S 6, we complete the proof of Theorem 1. One step is to show that
if $\tilde M$ is a regular $\bz_{p}$--cover which is not a simple
$\bz_{p}$--cover of $M$, then we may delete a
simple closed curve $\g$ in $M$ so that the inverse image of $M-\g$
is a simple $\bz_{p}$--cover of $M-\g$.
In \S 7, we  prove Theorem 2, Theorem 3,  Corollary 1, and Corollary 2.

 \section{Turaev--Viro invariants}

We also want to consider the associated Turaev--Viro invariants, which
one may define by $\text{\rm tv}_{r}(M)= w_{r}(M) \overline{w_{r}(M)}$.
Here conjugation is defined by the usual conjugation defined on the
complex numbers.
$\text{\rm tv}_{r}(M)$  was first defined as a state sum by
Turaev and Viro \cite{TV}, and later shown to be given by the above
formula
separately by Walker and Turaev.  A very nice proof of this fact was
given by Roberts \cite{R}. We will use a state sum definition in the
form used by Roberts.
We pick a triangulation of $M$, and  sum
certain contributions over $r$--admissible colorings $C$  of the triangulation.
A coloring of a triangulation assigns to each 1--simplex a nonnegative integral
color less than $r-1$. The coloring is admissible if:
for each 2--simplex the colors assigned to the three edges $a$, $b$, and $c$
satisfy $a +b+c$ is even, $a +b+c \le 2r-4$, and $|a-b| \le c \le a+b$.
$$
\text{\rm tv}_{r}(M)= \sum_{c \in C} \prod _{v\in V}\eta^{2} \
\prod _{e\in E}\Delta(c,e)\ \prod _{f\in F}\theta(c,f)^{-1}\
\prod _{t\in T}\ \text{Tet}(c,t)
$$
Here $V$ is the set of vertices,  $E$ is the set of edges (or 1--simplexes),
 $F$ is the set of faces (or 2--simplexes), and  $T$ is the set of
 tetrahedrons (or 3--simplexes) in the triangulation.
The contributions are products of certain evaluations in the sense of
Kaufman--Lins  \cite {KL} of colored planar graphs.
$\Delta(c,e)=\Delta_{c(e)}$, where  $\Delta_i$ is the evaluation of a loop
colored the color $i$.
$\theta(c,f)$ is the evaluation of an unknotted
theta curve whose edges are colored with the colors assigned  by $c$
to the
edges of $f$. $\text{Tet}(c,t)$ is the evaluation of a tetrahedron
whose edges are colored with the colors assigned  by $c$ to the
tetrahedron $t$.
However we take all these evaluations in $\bc$,  taking
$A^2$ to be $\l=\xi_{2r}$. Also $\eta^{2}= -\frac{(\l-\l^{-1})^{2}}{2r}$.
So $\text{\rm tv}_{r}(M)$ lies in
$\bz[\frac 1 {2r},\l] $.
This follows from the following lemmas
and the formulas in \cite{KL} for these evaluations.

\proclaim {Lemma 1} For $j$ not a multiple of $2r$,
 $(1-\l^{j})^{-1} \in \bz[\frac 1 {2r},\l]$.
\endproclaim
\demon  $\prod_{s=1}^{2r-1}(x-\l^{s})= \sum_{i=0}^{2r-1} x^{i}$.
Letting $x=1$, $\prod_{s=1}^{2r-1}(1-\l^{s})= 2r$. \qed \enddemo

\proclaim {Lemma 2} For $n \le r-1$, the ``quantum integers''
 $[n]= \frac{\l^{n}-\l^{-n}}
{\l-\l^{-1}}= \l^{1-n} \frac{\l^{2n}-1}
{\l^{2}-1}$ are units in $ \bz[\frac 1 {2r},\l]$.
\endproclaim

Since it is fixed by complex conjugation, $\text{\rm tv}_{r}(M) \in
\bz[\frac 1 {2r},\l+\l^{-1}]$.

\proclaim{Theorem 4} If the action is free, then
$$ \text{\rm tv}_r(\tilde M)\equiv (\text{\rm tv}_r(M))^{|G|} \pmod{p}.$$
\endproclaim

\demon  A chosen triangulation $T$ of $M$ lifts to a
triangulation $\tilde T$
of $\tilde M$. Each admissible coloring of $T$ lifts to an admissible
coloring of $\tilde T$. As each simplex of $M$ is covered by $|G|$
simplexes  of $\tilde M$, the contribution of a lifted coloring to
the sum for $\text{\rm tv}_r(\tilde M)$ is the $|G|$th power of the
contribution of the  original coloring to $\text{\rm tv}_r(M)$.
$|G|$ acts freely on the set of colorings of $\tilde T$ which are
not
lifts of some coloring of $T$. Moreover the contribution of each
such triangulation in a given orbit of this $G$ action is
constant.
Thus the contribution of the non-equivariant colorings is a multiple
of $p$. Making use of the equation $x^{p^s}+ y^{p^s}\equiv (x+y)^{p^s}
\pmod{p}$, the
result follows. \qed \enddemo

\proclaim{Lemma 3} If $r$ divides $\frac{p^s \pm1}2$, $\D_i\equiv (\D_i)^{p^s}
\pmod{p}$, for
all $i$, and
$\eta^2\equiv (\eta^2)^{p^s} \pmod{p}$.
\endproclaim

\demon  Since $r$ divides $\frac{{p^s} \pm1}2$, $2r$ divides ${p^s} \pm1$.
Thus
$\l^{{p^s} \pm 1}=1$, and $\l^{{p^s}}=\l^{\mp 1}$.
Thus
$$\D_{1}^{{p^s}}= \(-\l-\l^{-1}\)^{{p^s}}
\equiv -\l^{{p^s}}-\l^{-{p^s}}=
-\l-\l^{-1}=
\D_{1} \pmod{p}. $$
Also $\D_0=1$. Thus $\D_i\equiv (\D_i)^{p^s} \pmod{p}$, if $i$ is zero or
one. Using the recursion formula $\D_{i+1}= \D_1 \D_i- \D_{i-1}$,
$\D_i\equiv (\D_i)^{p^s} \pmod{p}$  follows by induction.
Here is the
inductive step:
$$\multline
(\D_{i+1})^{p^s}= \(\D_1 \D_i- \D_{i-1}\)^{p^s} \\
\equiv \D_1^{p^s} \D_i^{p^s}-
\D_{i-1}^{p^s}\equiv \D_1 \D_i- \D_{i-1}=\D_{i+1} \pmod{p}.\endmultline
$$
It follows that $\sum_{i=0}^{r-2} \D_i^2 \equiv \(\sum_{i=0}^{r-2}
\D_i^2\)^{p^s}
\pmod{p}$. As $\eta ^2\sum_{i=0}^{r-2} \D_i^2=1$, and
$\bz[\l,\frac 1 {2r}]\slash p(\bz[\l,\frac 1 {2r}])$ is a direct sum of
fields, we have $(\eta^2)^{p^s}\equiv \eta^2 \pmod{p}$.
\qed \enddemo

\proclaim{Theorem 5} If $r$ divides $\frac{p \pm1}2$,  then
$$ \text{\rm tv}_r(\tilde M)\equiv (\text{\rm tv}_r(M))^{|G|} \pmod{p}.$$
If $G$ is cyclic and acts semifreely and $r$ divides $\frac{{|G|} \pm1}2$,
the same congruence holds.
\endproclaim
\demon  We pick our triangulation of the base so that the image
of the fixed point set is a one dimensional subcomplex. By Lemma 3,
whether a colored  simplex in the base lies in the image of a simplex with a
smaller orbit
or not it contributes the same amount modulo $p$ to a product associated to
an equivariant coloring. Thus the proof of Theorem 4 still goes
through. \qed \enddemo

\section{Quantum invariants, $p_{1}$--structures, and signature\\
defects}

Let $M$ be a closed 3--manifold with a $p_1$--structure \cite{BHMV}.
A fat colored  graph in $M$ is a trivalent graph
embedded in $M$, with a specified 2--dimensional thickening (ie, banded in the
sense of
\cite{BHMV}) whose edges have been colored with nonnegative integers less than
$r-1$. At each vertex the colors on the edges $a$, $b$, and $c$ must satisfy
the  admissibility conditions:
$a +b+c$ is even, $ a +b+c \le 2r-4$, and $|a-b| \le c \le a+b$.
Let $J$ be such a graph (possibly empty) in $M$.
Recall the quantum invariant $\left<(M,J)\right>_{2r}\in k_{2r}$ defined in
\cite {BHMV}.
Consider the homomorphism
[\Cite{MR}; note page 134] $\psi\co k_{2r}\rightarrow \bc$ which sends $A$ to
$-\xi_{4r}$,
 and sends $\k$ to $\xi_8
\xi_{4r}^{-1}$.
Let $R_r$ denote the image of $\psi$. By abuse of notation let $\k$
denote $\psi(\k)$, and $\eta$ denote $\psi(\eta)$.  Let $t=4r$ if $r$
is even and $t= 8r$ if $r$ is
odd.  Then $R_r=\bz[ \frac 1 {2r}, \x_t]$.
Further abusing notation, we let $\left<(M,J)\right>$ denote $\psi
\left<(M,J)\right>_{2r}\in R_r$.  If $M$ is a
3--manifold without an assigned $p_1$--structure, we let $w_r(M,J)$
denote
$\left<(M',J)\right>$ where $M'$ denotes $M$ equipped
with a $p_1$--structure with $\sigma$--invariant zero.
If $M$ already is assigned a $p_1$--structure, we let $w_r(M,J)$ denote
$\left<(M',J)\right>$ where $M'$ denotes $M$ equipped with a reassigned
$p_1$--structure with $\sigma$--invariant zero.
One has that $w_r(M,J)= \k^{-\si(M)} \left<(M,J)\right>$. This agrees with
$w_r(M,J)$ as defined in the introduction.

Assume now that $M$ has been assigned a $p_1$--structure. Let $\tilde M$ be  a
regular $G$ covering space. Give
$\tilde M$ the induced $p_1$--structure, obtained by pulling back the
structure on $M$.
The following lemma  generalizes [\Cite{G2}; 3.5]. It does not require that $G$
be a $p$--group.

\proclaim{Lemma 4} $3 \  \text{\rm def }(\tilde M \rightarrow
M) = |G|\ \si(M)- \si(\tilde M)$. In particular
$3 \  \text{\rm def }(\tilde M \rightarrow
M)$ is an integer. \endproclaim

\demon  Pick a 4--manifold $W$ with boundary $|G|$ copies of $M$
such that the cover extends. We may connect sum on further copies of
$\bc P(2)$ or
$\overline{\bc P(2)}$ so that the $p_1$--structure on $M$ also
extends. Let $\tilde W$
be the associated cover of $W$ with boundary $|G|$ copies of $\tilde
M$.
We have
$$ |G|\ \si(M) = 3 \  \text{Sign}(W),\quad
 |G|\ \si(\tilde M) = 3 \  \text{Sign}(\tilde W),$$
$$|G|\ \text{\rm def }(\tilde M \rightarrow M)= |G|\ \text{Sign}(W)-
\text{Sign}(\tilde W).\eqno{\qed}$$\enddemo

Using this lemma, we  rewrite Theorem 1  in an equivalent form.
The conclusion is simpler. On the other hand, the hypothesis involves
the notion of a $p_1$--structure. Since $p_1$--structures are sometimes
a stumbling block to novices, we stated our results in the
introduction
without reference to $p_1$--structures.

\proclaim{Theorem 1$^\prime$} Let $M$ have a $p_1$--structure, and
$\tilde
M$ be a regular $G$ cover of $M$ with the induced
$p_1$--structure. Then
$$\left<\(\tilde M,\tilde J\)\right>\equiv  \left<(M,J)\right>^{|G|}
\pmod{p}.$$ \endproclaim

Note that we may define $\left<(M,\Cal J)\right>$ for $\Cal J$ a linear
combination over $R_{r}$ of fat colored graphs in $M$, by extending
 the function
$\left<(M,J)\right>$ linearly. If $\Cal J= \sum_i a_i J_i$, we define $\tilde
\Cal J$ to be $\sum_i a_i^{|G|} \tilde J_i$.
Since the pth power map is an automorphism of
$R_{r}
\slash pR_{r}$,
we have that if Theorem 1$^{\prime}$ is true for a given type
manifold $M$,
then it is true for such manifolds when we replace $J$ and $\tilde J$
by linear combinations over $R_r$ of colored fat graphs:
$\Cal J$ and $\tilde {\Cal J}$.

Finally we note that if Theorem 1$^{\prime}$ is true for $G= \bz_p$, then it
will follow
for $G$ a general finite $p$--group. In the next three sections we prove it
for $G= \bz_p$.

\section{Simple unbranched $\bz_{p^s}$--covers}

A regular $Z_{p^s}$--covering space $\tilde X$ of $X$ is classified by an
epimorphism $\phi\co 
H_1(X)$\break $ \rightarrow \bz_{p^s}$. If $\phi$ factors through $\bz$, we say
$\tilde X$ is a simple $\bz_{p^s}$--cover. In this section, we prove
Theorem 1$^{\prime}$ for simple $\bz_{p^s}$--covers. We also obtain a version
for
simple $\bz_{p^s}$--covers of manifolds whose boundary is a torus.

If $\chi\co  H_1(M) \rightarrow \bz$ is an epimorphism, let $\chi_{p^s}\co 
H_1(M) \rightarrow \bz_{p^s}$
denote the composition with reduction modulo ${p^s}$.
Suppose $\tilde M$ is classified by $\chi_{p^s}$.
Consider a Seifert surface for $\chi$ ie, a closed surface  in $M$
which is Poincare dual to $\chi$. We may and do assume that this
surface is in general position with respect to the colored fat graph
$J$. Then the intersection of $J$ with $F$ defines some banded
colored points.
This surface also acquires a
$p_1$--structure. Thus $F$ is an object in the cobordism category
$\Cal
C^{p_1,c}_{2,r-1}$, [\Cite{BHMV}; 4.6].
Let $E$ be the cobordism from $F$ to $F$ obtained by slitting $M$
along $F$.
We view $E$ as a morphism from $F$ to $F$ in the cobordism category
$\Cal C^{p_1,c}_{2,r-1}$.  Then $M$ is the mapping torus of $E$,
and $\tilde M$ is the mapping torus for $E^{p^s}$.

We may consider the TQFT  which is a functor  from $\Cal
C^{p_1,c}_{2,r-1}$, to the category of modules over $R_r$ obtained
taking $p=2r$ in \cite{BHMV} and applying the change of coefficients
$\psi\co k_{2r}\rightarrow R_r$.

By [\Cite{BHMV}; 1.2], we have
$$\left<(M,J)\right> =\text{Trace}\( Z(E) \)\quad \text{and}\quad
\left<\(\tilde M,
\tilde J\)\right>
=\text{Trace} \(Z\(E^{p^s}\)\).$$
Let $\Cal E $ be the matrix for $Z(E)$ with respect to some basis for
$V(F)$. $\Cal E $ has entries in $R_{r}=\bz[\frac 1 {2r}, \xi_{t}]$.
Write the entries as polynomials in $\xi_{t}$ whose coefficients are
quotients of integers by powers of  $2r$. Let $v \in \bz$ such
that $2r v\equiv 1 \pmod{p}$. Let $\Cal E' $ denote the matrix over
$\bz[\xi_{t}]$ obtained by replacing all powers of $2r$ in the
denominators of entries by powers of $v$ in the numerators of
these
entries. We have
$$\text{Trace}\( Z(E) \)=\text{Trace}\(\Cal E \) \equiv
\text{Trace}\(\Cal E' \) \pmod{p} \text{ and},$$
$$ \text{Trace}\( Z(E)^{p^s} \)=
 \text{Trace}\( \Cal E^{{p^s}} \) \equiv \text{Trace}\(\Cal  E^{\prime {p^s}}
\) \pmod{p}.$$
However all the eigenvalues of $\Cal E^\prime$ are themselves
algebraic
integers. The trace of  $\Cal E^{\prime}$ is the sum of these
eigenvalues
counted with multiplicity. The trace of  $\Cal E^{\prime {p^s}}$ is the
sum of
${p^s}$th powers of these eigenvalues counted with multiplicity.
Therefore
$$\text{Trace}\( \Cal E^{\prime {p^s}} \) \equiv
\(\text{Trace}\( \Cal E^{\prime} \)
\)^{{p^s}} \pmod{p}.$$ Putting these equations together proves Theorem
1$^\prime$ for simple $\Bbb Z_{p^s}$ covers.

We now wish to obtain a version of Theorem 1$^\prime$ for manifolds
whose
boundary is a torus. Let $N$ be a compact oriented 3--manifold with
$p_1$--structure with boundary
$S^{1} \times S^{1}$.  $S^{1} \times S^{1}$ acquires a
$p_1$--structure  as the boundary. Let $J$ be a colored fat  graph in
$N$ which is disjoint from the boundary. Then $(N,J)$ defines an
element of $V(S^{1} \times S^{1})$ under the above TQFT. We denote
this element by $[N,J]$.

If $\chi\co  H_1(N) \rightarrow \bz$, let $\chi_{p^s}\co  H_1(N) \rightarrow
\bz$ denote the composition with reduction modulo ${p^s}\co\bz
\rightarrow \bz_{p^s}$. Suppose $\chi_{p^s}$ restricted to the boundary is
an epimorphism.  Suppose $\tilde N$ is a regular $Z_{p^s}$ covering
space given by $\chi_{p^s}$. $\tilde N$ has an induced
$p_1$--structure. Suppose  that we have identified the boundary with
$S^{1} \times S^{1}$ so that $\chi$ restricted to the boundary is
$\pi_{1*}\co  H_1(S^{1} \times S^{1}) \rightarrow H_1(S^1)$, followed by
the standard isomorphism. Here $\pi_1$ denotes projection on the
first factor.
We can always identify the boundary in this way.

Let $\tilde N$ be the regular $Z_{p^s}$ covering space given by
$\chi_{p^s}$, and $\tilde J$ the colored fat graph in $\tilde N$ given by
the
inverse image of $J$. The boundary of $\tilde N$ is naturally
identified with
$\tilde {S^1} \times S^1$. Equip  $S^1 \times D^2$ with a
$p_1$--structure
extending the $p_1$--structure that $S^1 \times S^1$ acquires as the
boundary of $N$.
Equip  $\tilde {(S^1)} \times D^2$ with a $p_1$--structure
extending the $p_1$--structure $\tilde {(S^1)} \times S^1$ acquires as
the
boundary of $\tilde N$. This is the same $p_1$--structure it gets as
the cover of $S^1 \times S^1$. We have
$[\tilde N,\tilde J]\in V(\tilde {S^1}\times S^1 )$. Also  $\partial (S^1\times
D^2 )= \partial N$, and
$\partial (\tilde {S^1}\times D^2)= \partial \tilde N$.
For $0 \le i \le r-2$, let $e_i =[S_i]\in V(S^1 \times S^1)$, where
$S_i$ is
 $S^1 \times D^2$ with the core with the standard thickening  colored $i$.
Similarly let $\tilde e_i=[\tilde S_i]\ \in V( \tilde{S^1} \times
S^1)$,
where $\tilde S_i$ is $\tilde S^1 \times D^2$ with the core with the
standard
thickening colored $i$.

\proclaim{Proposition 1} If $[N, J]= \sum_{i=0}^{r-2} a_i e_i, $ and
$[\tilde N, \tilde J]= \sum_{i=0}^{r-2} \tilde a_i \tilde e_i, $
 then $\tilde a_j\equiv  a_j^{p^s} \pmod{p}$, for  all $j, $ such that $0
\le j \le r-2$.
\endproclaim

\demon  Note
$(\tilde N,\tilde J) \cup_
{ \tilde {S^1}\times S^1 }
-\tilde {S_j}$
is a simple cover of $( N, J )\cup_{{S^1}\times S^1} -S_j$.
Moreover $\left<(\tilde N,\tilde J) \cup_{\tilde {S^1}\times S^1}- \tilde
S_j\right>= \tilde a_j$
and $\left< (N, J) \cup_{ {S^1}\times S^1} -S_j\right>=a_j$. By Theorem
1$^\prime$,
$\tilde a_j\equiv  a_j^{p^s} \pmod{p}$.
\qed
\enddemo

\proclaim{Proposition 2} If $\Cal J$ is a linear combination over
$R_{r}$ of colored fat graphs in $S^{1} \times D^{2}$  then $[S^{1} \times
D^{2},\Cal J] \in
V(S^{1} \times S^{1})$ determines
 $[\tilde S^{1} \times D^{2},\tilde \Cal J] \in V(\tilde S^{1} \times
S^{1})$
 modulo $p$.
\endproclaim
\demon  Just take  $N=S^{1} \times D^{2}$, and sum over the
terms in $\Cal J$.
\qed
\enddemo

See the paragraph following the statement of Theorem 1$^\prime$ for the
definition   of $\tilde \Cal J$.

\section{$\bz_p$--covers of lens spaces}

$L(m,q)$ can be described as $-m \slash q$ surgery to an unknot in $S^{3}$.
A meridian of this unknot becomes a curve in $L(m,q)$ which we refer
to as a meridian of $L(m,q)$.
Below, we verify directly that Theorem 1 (and thus Theorem
1$^\prime$) holds
 when $M$ is a lens space, $J$ a meridian colored $c$, and $G=\bz_p$.
  By general
 position any fat graph $J$  in a lens space can be isotoped into a
tubular neighborhood
 of any meridian. Without changing the invariant of the lens space
 with $J$ or the cover of the lens space with $\tilde J$ one can
 replace $J$ by a linear combination of this meridian with various
colorings and $\tilde
 J$ by the same linear combination of the inverse image of this
  meridian with the same
colorings in the covering space. This
follows from Proposition 2. Thus it will follow that Theorem 1  (and
thus Theorem 1$^\prime$) will hold if  $M$ is a lens space, and $J$
is  any fat colored graph in $M$, and $G=\bz_p$.
This is a step in the proof of Theorem 1$^\prime$ for $G=\bz_p$.

Consider the $p$--fold cyclic cover $L(m,q) \rightarrow L(mp,q)$.  We
assume $m$, $q$ are  greater than zero.  $q$ must be relatively prime
to $m$ and $p$.

\proclaim{Lemma 5}
  $$\text{\rm def }\(L(m,q)   \rightarrow L(mp,q)\)= \frac 1 m \(
\text{\rm def }\(S^3   \rightarrow L(mp,q)\)- \text{\rm def }\(S^3
\rightarrow L(m,q)\)\).$$
\endproclaim

\demon Suppose $Z \rightarrow X$ is a regular $\bz_{mp}$
covering of 4--manifolds with boundary, and on the boundary we have
$mp$ copies of the regular $\bz_{mp}$ covering $S^{3} \rightarrow
L(mp,q)$. Let $Y$ denote  $Z$ modulo the action of $\bz_{m}\subset
\bz_{mp}$.
Then  $Z \rightarrow Y$ is a regular $\bz_{m}$
covering, and on the boundary we have
$mp$ copies of the regular $\bz_{m}$ covering $S^{3} \rightarrow
L(m,q)$. Moreover $Y \rightarrow X$ is a regular $\bz_{p}$ covering
and on the boundary we have
$mp$ copies of the regular $\bz_{p}$ covering $L(m,q) \rightarrow
L(mp,q)$. We have:
$$mp\(\text{\rm def }(L(m,q)   \rightarrow L(mp,q)\)= p\ \text{Sign}(X)-
\text{Sign}(Y),$$
$$mp\(\text{\rm def }(S^{3}   \rightarrow L(m,q)\)= m\ \text{Sign}(Y)
-\text{Sign}(Z)\quad  \text{and}$$
$$mp\ \text{\rm def }(S^{3}   \rightarrow L(mp,q)= mp\ \text{Sign}(X)-
\text{Sign}(Z).$$
The result follows.
\qed
\enddemo

Suppose $\Cal H$ is a normal subgroup of a finite group  $\Cal G$ which acts
freely on $M$.  By the above argument, one can show more generally  that:
$$ \text{\rm def }\(M   \rightarrow M/\Cal G\)=
\text{\rm def }\(M   \rightarrow M/\Cal H\)+ |\Cal H|\ \text{\rm def }\(M/\Cal
H   \rightarrow M/\Cal G\).
$$
\noindent According to Hirzebruch \cite{H} (with different conventions)
$$3\ \text{\rm def }\(S^3   \rightarrow L(m,q)\)= 12 m\ s(q,m)\in \bz.$$  See
also
[\Cite{KM2}; 3.3], whose conventions we follow.

Thus
$$3\ \text{\rm def }\(L(m,q)   \rightarrow L(mp,q)\)= \frac 1 m \(
12 m p\ s(q,mp)- 12 m\ s(q,m)\) \in \bz.$$
\noindent We will need  reciprocity for generalized Gauss sums in the form due
to Siegel [\Cite{BE}; Formula 2.8]. Here is a slightly less general
form which suffices for our purposes: $$ \sum_{k=0}^{ \g
-1} \xi^{\a k^2 +\b k}_{2 \g}= \xi_8 \ \xi^{-\b^2}_{8\a\g} \sqrt{
\frac \g \a } \ \
\sum_{k=0}^{ \a -1} \xi^{-(\g k^2 +\b k)}_{2 \a}$$
\noindent where $\a,\b,\g \in \bz$ and $\a$, $\g > 0$ and $\a \g+ \b$
is even.

\noindent We use this reciprocity to rewrite the  sum $$
\align
\sum_{n=1}^{m}\xi_m^{(qr)n^2+(q l \pm 1)n}
&= \sum_{n=1}^{m}\xi_{2m}^{(2qr)n^2+2(q l \pm 1)n}\\
&= \xi_8 \ \xi_{16mqr}^{- 4 (ql \pm 1)^2} \sqrt{ \frac m {2 qr}}
\sum_{n=1}^{2qr}\xi_{4qr}^{-mn^2-2(q l \pm 1)n}.\\
\endalign $$
\noindent We substitute this into a formula from [\Cite{G1}; section 2].
$$
w_r( L(m,q),\mu_c) = \frac { i (-1)^{c+1}} {\sqrt{2rm}} \sum_{\pm}\pm
\xi_{4rmq}^{mb - mq \Phi(U) + q^2 l^2 +1 \pm 2q l}
\sum_{n=1}^{m}\xi_m^{(qr)n^2+(q l \pm 1)n},$$
\noindent where $U=(\smallmatrix q & b \\ m&
d\endsmallmatrix)$, and $l=c+1$.  We remark
that the derivation given in \cite{G1} is valid for $m \ge 1$.
Also using
$$mb - mq \Phi(U) + q^2 l^2 +1=q^2(l^2-1)+12 mq\ s(q,m)$$
\noindent again from \cite {G1}, we obtain:
$$
w_r( L(m,q),\mu_c) = \frac { \xi_8^3 (-1)^{c+1}} {2r\sqrt{q}}
\xi_{4rm}^{12 m\  s(q,m)} \xi_{4rmq}^{-q^2-1}
\sum_{\pm}\pm \sum_{n=1}^{2qr}\xi_{4qr}^{-mn^2-2(q l \pm 1)n}.$$
\noindent Similarly
$$
w_r( L(mp,q),\mu_c) = \frac { \xi_8^3 (-1)^{c+1}} {2r\sqrt{q}}
\xi_{4rmp}^{12 mp\  s(q,mp)}\xi_{4rmpq}^{-q^2-1}
\sum_{\pm}\pm \sum_{n=1}^{2qr}\xi_{4qr}^{-mpn^2-2(q l \pm 1)n}.
$$
\noindent So we now work with congruences modulo $p$  in the ring of algebraic
integers after we have inverted $2r$ and $q$.
We need to show that
$$ w_r( L(m,q),\mu_c)\equiv
\k^{3\  \text{\rm def }\(L(m,q)   \rightarrow L(mp,q)\)}
 \(w_r( L(mp,q),\mu_c)\)^p \pmod{p}. \eqno{(5.1)}$$
\noindent We have that:
$$\align \k^{3\  \text{\rm def }\(L(m,q)   \rightarrow L(mp,q)\)}&=
\(\xi_8 \xi_{4r}^{-1}\)^{\frac 1 m \(
12 mp\  s(q,mp)- 12 m\ s(q,m)\)}\\
 &=\xi_8^{\frac 1 m \(
12 mp\  s(q,mp)- 12 m\ s(q,m)\)}
\xi_{4rm}^{12 m\  s(q,m)- 12 mp\ s(q,mp)}, \endalign $$
$$\(\xi_{4rmpq}^{-q^2-1}\)^p =\xi_{4rmq}^{-q^2-1}$$
$$ (-1)^p = -1 \text{\ \ and},$$
$$(2r)^p \equiv 2r \pmod{p}.$$
\noindent Note that
$\xi_{4qr}^{-mn^2-2(q l \pm 1)n}$ only depends on $n$ modulo $2qr$.

\noindent Thus $$\align \(\sum_{n=1}^{2qr}\xi_{4qr}^{-mpn^2-2(q l \pm 1)n}\)^p
& \equiv \sum_{n=1}^{2qr}\xi_{4qr}^{-m(pn)^2-2(q l \pm 1)pn}\\
&\equiv \sum_{n=1}^{2qr}\xi_{4qr}^{-mn^2-2(q l \pm 1)n} \pmod {p} .
\endalign 
$$
\noindent Here we have made use of the fact that as $n$ ranges over
all the congruence classes modulo $2qr$ so does $pn$. Let $\binom q
p$ denote the Legendre--Jacobi symbol.
$$(\sqrt q)^p = q^{\frac {p-1} 2}(\sqrt q)\equiv
\binom q p \sqrt q \pmod{p}$$
 \noindent Equation $(5.1)$ will follow when we show:
$$\binom q p \xi_8^3 =
\xi_8^{\frac 1 m \(12 mp\ s(q,mp)-12 m\ s(q,m)\)}\xi_8^{3p},$$
or
$$ \xi_8^{3-3p
+\frac 1 m \(12 m\ s(q,m)-12 mp\ s(q,mp)\)}=\binom q p .
\eqno{(5.2)}$$
\noindent We will make use of a congruence of Dedekind's [\Cite{R}; page 160
(73.8)]:
for $k$ positive and odd:
$$12k\ s(q,k) \equiv k+1 - 2 \binom q k \pmod{8}.$$
\noindent We first consider the case that  $m$ is odd. We have:
 $$\align 12 m\ s(q,m)-12 mp\ s(q,mp)&\equiv  m (1-p) + 2 \(\binom q
{mp}-\binom q m\)\\
&=  m (1-p) + 2 \binom q m \(\binom q {p}-1\)
\pmod{8}. \endalign $$
\noindent If $m$ is odd, the equation  (5.2)  becomes:
$$(3-3p)m+ m (1-p) + 2 \binom q m \(\binom q {p}-1\)\equiv
\cases
0,&\text{for $\binom q p =1$}\\
4,&\text{for $\binom q p =-1$}
\endcases
\pmod{8},$$
\noindent which is easily checked.
 For the rest of this section, we consider the case that  $m$ is even.
In this case, $q$ must be odd. First we use Dedekind reciprocity
[\Cite{R}; page 148 (69.6)] to rewrite
$  \frac 1 m \(12 m\ s(q,m)-12 mp\ s(q,mp)\)$ as:
$$ \multline \frac 1 {mq} (-12mq\ s(m,q) +m^2 +q^2+1-3mq+12mqp
\ s(mp,q)\\ -m^2p^2-q^2-1+3mpq )\endmultline$$
$$  \text{or \ } \frac 1 q \(-12q\ s(m,q) +m -3q+12qp\ s(mp,q)-m p^2+3pq
\).$$
As $q$ is odd, $q^2= 1 \pmod{8}$.  So modulo eight the above
expression is
$$q \(-12q\ s(m,q) +m -3q+12qp\ s(mp,q)-m p^2+3pq \).$$
As $p$ is odd, $p^2= 1 \pmod{8}$, and $m-mp^2=0 \pmod{8}$.
The expression, modulo eight, becomes
$$-12q^2\ s(m,q)  -3q^2+12q^2 p\ s(mp,q)+3pq^2 . $$
So the exponent of $\x_{8}$ in Equation (5.2)
modulo eight is:
$$q\( -12q\ s(m,q) +12q p\ s(mp,q) \). \eqno{(5.3)}$$
\noindent Again using Dedekind's congruence, we have:
$$12 q \ s(m,q) \equiv q+1 - 2 q \binom m q \pmod{8}.$$
\noindent Similarly
$$12 q \ s(mp,q) \equiv q+1 - 2 q \binom {mp} q \pmod{8}.$$
\noindent The  expression (5.3)  modulo eight, becomes
 $$q(p-1)(q+1) +2 \binom m q \(1- p \binom p q \)$$
\noindent Thus we only need to see
 $$\xi_8^{q(p-1)(q+1) +2 \binom m q (1- p \binom p q )}=\binom q
p.$$
Using quadratic reciprocity for Jacobi--Legendre symbols, this
becomes
$$(-1)^
{\frac {p-1} 2 \frac {q+1} 2 +   \binom m q \frac {\(1-p \binom p q
\)}2 } = (-1)^{\frac {p-1} 2 \frac{q-1} 2}\binom p q,$$
or $$(-1)^
{\frac {p-1} 2  + \frac {\(1-p \binom p q \)}2 } = \binom p q,$$
which is easily checked.

\section{Unbranched $\bz_p$--covers}

\proclaim{Lemma 6} If $\phi \in H^1(M, \bz_p)$ is not the
reduction of an integral class, then there is a simple closed curve
$\g$  in $M$ such that
the  restriction of
$\phi $ to $M-\g$ is the reduction of an integral class.\endproclaim

\demon  We consider the Bockstein homomorphism $\b$ associated
to the short exact sequence of coefficients:
$0 \rightarrow \bz \rightarrow \bz\rightarrow \bz_p \rightarrow 0$.
Let $\g$ be a simple closed curve which represents the element which
is Poincare dual
to $\b(\phi) \in H^2(M,\bz)$. The restriction of $\b(\phi)$ to
$H^2(M-\g,\bz)$
is zero. This is easily seen using the geometric description of
Poincare duality which makes use of dual cell decompositions.

Naturality of the long exact Bockstein sequence completes the proof.
\qed\enddemo

It follows that any non-simple unbranched cyclic covering of a closed
3--manifold can be decomposed as the union of two simple coverings one
of which is a covering of a solid torus. Let $N$ denote $M$ with a
tubular neighborhood of $\g$ deleted. We may replace $N$ in $M$ by a
solid torus with a linear combination of colored cores obtaining a new
manifold $(M',J')$, such that $\left<(M',J')\right>=\left<(M,J)\right>$. By
Proposition 2,
 $\left<\(\tilde M',\tilde J'\)\right>=\left<\(\tilde M, \tilde J\)\right>$.
$M'$ is a union
 of two solid tori, and so is a
lens space or $S^{1} \times S^{2}$. Note that any cover of $S^{1}
\times S^{2}$ is simple.  As Theorem 1$^\prime$ when $G=\bz_p$ has already
been established for $M'$,  we have now
established Theorem 1$^\prime$
when $G=\bz_p$.

\section{Branched $\bz_{p^s}$--covers }

We need to refine the last statement of Lemma 3:

\proclaim{Lemma 7} If $r$ divides
$ \frac{{p^s}\pm1}2$, $\eta^{{p^s}} \equiv  \mp {\binom {-2r} p}^s \eta
\pmod{p}.$
\endproclaim

\demon  We have that
$\eta= -\frac {(\l-\l^{-1})}{\sqrt{2r}}i$,
$\l^{{p^s}} \equiv \l^{\mp 1} \pmod{p}$,
$i^{{p^s}}= {\binom {-1} p}^s i$,
$(\sqrt{2r})^{{p^s}}\equiv  {\binom {2r} p}^s \sqrt{2r}\pmod{p}$.
\qed
\enddemo

We  need the following  which follows from Proposition 1,
Lemma 3, and Lemma 7. In the notation developed at the end of \S 4 ,
let $\o= \eta \sum_{i=o}^{r-2} \D_{i} S_{i}$.

\proclaim{Proposition 3} If $r$ divides
$ \frac{p^s \pm1}2$, $\tilde \o \equiv \mp {\binom {-2r} p}^s \eta
\sum_{i=o}^{r-2}
\D_{i}
\tilde S_{i} \pmod{p}$.
\endproclaim

\medskip {\bf 7.1\qua Proof  of Theorem 2 for $G$ cyclic} \medskip

Suppose $M$ is a 3--manifold, $L$ is a link in this 3--manifold
and
there exists a homomorphism $\phi\co H_{1}(M - L) \rightarrow \bz_{p^s}$
which
sends each meridian of $L$ to a unit of $\bz_{p^s}$. Then we may form a
branched cover  $\tilde M$ of $M$ branched along $L$. Every
 semi-free $\bz_{p^s}$
action on an oriented manifold arises in this way. Then we may pick
some
parallel curve to each component of $L$ whose homology class maps to
zero. Perform
integral surgery to $M$ along $L$ with framing given by these
parallel curves, to form $P$. Then we may complete the
regular
unbranched cover of $M-L$ given by $\phi$ to a  regular
unbranched cover $\tilde P$ of $P$. If we then do surgery to $P$
along an original
meridian (with the framing given by a parallel meridian) of each
component of $L$,
we recover $M$. Similarly if we do surgery to $\tilde P$ along
the inverse images of these meridians of $L$ then we recover $\tilde
M$.
Note that the inverse image of each meridian of $L$ is a single
component in $\tilde M$ and $\tilde P$. We give $M$  a
$p_{1}$--structure with $\si$ invariant zero.
Then  $P$ receives a $p_{1}$--structure as the result of $p_{1}$--surgery on $M$
[\Cite{BHMV}; page 925]. $\tilde P$ receives  $p_{1}$--structure as the cover of
$P$.  $\tilde M$ receives a $p_{1}$--structure as the result of $p_{1}$--surgery
on $\tilde P$.

Now let $J$ denote colored fat graph in $M$
disjoint from $L$. Now let $\Cal J^{+}$ denote the linear
combination of colored fat graphs in  $P$ given by $J$ together the
result of
replacing the meridians of $L$ by $\o$. As usual in this subject, the
union of  linear combinations is taken to be the linear combination
obtained by expanding multilinearly. Then
$\left<\(M,J\)\right>=\left<\(P,\Cal J^{+}\)\right>$, by \cite{BHMV}. By
Theorem 1,
$\left<\(\tilde P,
\tilde {\Cal J^+}\)\right> \equiv \left<\(P,{\Cal J^{+}}\)\right>^{p^s}
\pmod{p}$.
Using Proposition 3 and by \cite{BHMV},
we have
$\left< \(\tilde M,\tilde J\)\right>=  \(\mp {\binom {-2r}{p}}^s\)^{\b_0(L)}
\left<\(\tilde P, \tilde {\Cal J^+}\)\right>$.
As some power of $\k^{3}$ is minus one, and changing the
$p_1$--structure on $\tilde M$, has the effect of
multiplying $\left<(M,J)\right>$ by a power of $\k$,
this yields Theorem 2 for semifree actions of cyclic groups.
\qed

\medskip {\bf 7.2\qua  Proof of Theorem 2 for general $p$--groups.} \medskip
Now we assume $r$ divides $\frac {p \pm 1}2$. Thus we have the congruence for
every $\bz_p$ action by 7.1. However we can write the projection
from $\tilde M$ to $M$ as a sequence of quotients of $\bz_p$ actions.
\qed

\medskip {\bf 7.3\qua  Proof of Theorem 3 } \medskip

In the  argument of 7.1, if $M$ is a homology sphere, and $L$ is a
knot $K$, then the longitude of $K$ maps to zero under $\phi$, and
$P$ is obtained by zero framed surgery along $L$.  $\tilde M$ is a
rational homology sphere, and $\tilde P$ is obtained by zero framed
surgery along the lift of $K$.
The trace of both surgeries have signature zero.
If  we give $M$ a
$p_{1}$--structure with $\si$ invariant zero, then  $P$ also has a
$p_{1}$--structure with $\si$ invariant zero. Also $\tilde M$ and
$\tilde P$ have $p_{1}$--structures with the same  $\si$ invariant.
$\tilde P$ has  a $p_{1}$--structure  with $\si$ invariant $3\
\text{\rm def }(\tilde P \rightarrow P)$ by Lemma 4, but this is $-3
\si_{p^s}(K)$. Thus $ \(w_{r}(M,J)\)^{p^s} =
\left<\(M,J\)\right>^{p^s} \equiv
\mp {\binom {-2r}{p}}^s\left<\(\tilde M,
\tilde J\)\right> =
\mp (\binom {-2r}{p})^s \k^{3 \si_{p^s}(K)} w_{r}(\tilde M,\tilde J)
$,
\noindent modulo $p$.  This proves Theorem 3. \qed

\medskip {\bf 7.4\qua  Proof of Corollary 1 } \medskip

We obtain Corollary 1 from Theorem 3 by taking $M$ to be $S^3$,  $K$ to be the
unknot and $J$ to be $L$ colored one with the  framing given by a Seifert
surface for $J$. Using [\Cite{KL}; page 7], one has that
$$w_{r}(M,J) = \eta \  \left<L\right>_{A= -\xi_{4r}}=
\eta \  V_L\(\xi_{r}^{-1}\)\(-\xi_{2r}-\xi_{2r}^{-1}\).$$
In evaluating the Jones polynomial at $\xi_{r}^{-1}$, we  choose
$\sqrt{\xi_{r}^{-1}}$ to be $\xi_{2r}^{-1}$ (for the time being.)
Since the induced framing of $\tilde L$ is also given by the Seifert surface
which is the inverse image of the Seifert surface for $L$,
$$w_{r}\(\tilde M,\tilde J\) = \eta \  \left<\tilde L\right>_{A= -\xi_{4r}}=
\eta \  V_{\tilde L}\(\xi_{r}^{-1}\)\(-\xi_{2r}-\xi_{2r}^{-1}\).$$
 By Theorem 3, Lemma 3 and Lemma 7,
$V_{\tilde L}(\xi_{r}^{-1})\equiv \(V_L(\xi_{r}^{-1})\)^{p^s} \pmod{p}$.  This
means that the difference is $p$ times an algebraic integer.
$V_L(t)$ is a Laurent polynomial in $\sqrt{t}$ with integer coefficients.  We
have that
$$\(\sqrt{\xi_{r}^{-1}}\)^{{p^s}}=
\xi_{2r}^{-{p^s}}=
 \xi_{2r}^{\pm 1}=
\sqrt{\xi_{r}^{\pm 1}}
$$
Thus:
$$\(V_L(\xi_{r}^{-1})\)^{p^s}  \equiv   V_L(\xi_{r}^{\pm 1})\pmod{p}.$$
and so
$$V_{\tilde L}(\xi_{r}^{-1})\equiv  V_L(\xi_{r}^{\pm 1})\pmod{p}.$$
As all primitive $2r$th roots of unity are conjugate  over $\Bbb Z$,
$$V_{\tilde L}(\z^{-1})\equiv  V_L(\z^{\pm 1})\pmod{p}$$ holds if $\z$ is any
primitive rth root such that $r$ divides $\frac {{p^s} \pm 1} 2$, and $r > 2$.
We must choose the same $\sqrt{\z}$ when evaluating both sides.

For any link $\Cal L$, let $\#(L)$
 denote the number of components of $\Cal L$. One has that
$V_{\Cal L}(1) = (-2)^{\#(L)-1}$. One has that  $\#(\tilde L)\equiv \#(L)
\pmod{p-1}$.
So $V_{\tilde L}(1)  \equiv   V_L(1) \pmod{p}$.
Thus the stated congruence
holds if  $\z$ is any $\frac {{p^s} \pm 1} 2$-th root of unity, and $\z \ne
-1$.
\qed

\medskip {\bf 7.5\qua  Proof of Corollary 2 } \medskip
It suffices to prove the result for $n=p^s$ where $p^s \equiv \pm 1\pmod{8}$.
We apply Corollary 1 with $\z=i$ and choose
$\sqrt{i} = e^{\frac {10 \pi i} 8}$. By H. Murakami \cite{Mu}, for any proper
link
$V_{\Cal L}(i) = (-1)^{\text{\rm Arf}(\Cal L)}(\sqrt{2})^{\#(L)-1}$, with the
above choice of $\sqrt{i}$. Since $2$ is a square modulo $p$, $\sqrt{2}^{p-1}
=1$. Thus $(\sqrt{2})^{\#(L)-1}\equiv (\sqrt{2})^{\#(\tilde L)-1}\pmod{p}$.
Thus we conclude  $(-1)^{\text{\rm Arf}( L)}=(-1)^{\text{\rm Arf}(\tilde  L)}
\pmod{p}$.
Therefore $\text{\rm Arf}( L)\equiv \text{\rm Arf}(\tilde  L)\pmod{2}$.
\qed

\Addresses\recd


\Refs


\ref \key {BE} \by B Bernd\by R Evans \paper The determination of
Gauss Sums \jour Bulletin AMS,  \vol 5 \pages 107--129 \yr 1981\endref


\ref \key {BHMV} \by C Blanchet\by N Habegger\by G Masbaum\by  P Vogel
\paper
Topological quantum field theories derived from the Kauffman bracket
\jour
Topology,  \vol 34  \yr 1995 \pages 883--927\endref


\ref \key {G1} \by P Gilmer \paper Skein theory and
Witten--Reshetikhin--Turaev
Invariants of links in lens spaces \jour 
Commun. Math. Phys. 202 (1999) 411--419\endref

\ref \key {G2} \by P Gilmer\paper Invariants for one-dimensional cohomology
classes arising from TQFT \jour Topology and its Appl. \vol 75 \yr 1997
\pages 217--259 \endref



\ref \key {H} \by F Hirzebruch \paper The signature theorem:
reminiscences and
recreation \inbook Prospects of Math \moreref Ann. of Math.  Stud. Vol.
70, \publ
Princeton Univ. Press \yr 1971 \pages 3--31\endref


\ref \key {KL} \by L Kauffman\by S Lins \book Temperely Lieb recoupling
theory and invariants of 3--manifolds\moreref
Ann. of Math. Stud. Vol. 134,
 \publ Princeton Univ. Press \yr 1994\endref

\ref \key {KM1} \by R Kirby\by P Melvin \paper The 3--manifold
invariants of Reshetikhin--Turaev for SL(2,$\bc$)\jour Invent.
Math.\vol 105 \yr 1991\pages 473--545 \endref



\ref \key {KM2} \by R Kirby\by P Melvin \paper Dedekind sums,
$\mu$--invariants, and the signature cocycle \jour Math. Ann.\vol
299 \yr 1994 \pages 231--267 \endref

\ref \key {L} \by R Lickorish \paper The skein method for
three-manifold invariants \jour Knot Theory and its Ram. \vol 2 \yr 1993
\pages 171--194
\endref



\ref \key {MR} \by G Masbaum\by  J Roberts  \paper On central
extensions of
mapping class groups \jour Math. Ann. \vol 302 \yr 1995 \pages
131--150
\endref

\ref \key {M} \by K Murasugi \paper Jones polynomials of periodic links \jour
Pacific J. Math. \vol 131 \yr 1988 \pages
319--329
\endref

\ref \key {Mu} \by H Murakami \paper A recursive calculation of the Arf
invariant of a link \jour  J. Math. Soc. Japan, \vol 38 \yr 1986 \pages
335--338
\endref



\ref \key {RT}  \by N Reshetikhin\by  V Turaev \paper
Invariants of
3--manifolds via link polynomials  and quantum groups \jour Invent.
Math.\vol 103 \yr 1991 
\pages 547--597  \endref


\ref \key {R} \by H Rademacher \book Topics in analytic number theory
\publ Springer--Verlag \yr 1973  \endref



\ref \key {T} \by P Traczyk \paper $10_{101}$ has no period 7: a
criterion for periodic links\jour Proc. AMS,\vol 108 \yr
1990\pages 845--846
\endref

\ref \key {TV} \by V Turaev\by O Viro \paper State sum invariants of
3--manifolds and quantum 6j--symbols \jour Topology, \vol 31 \yr
1992\pages 865--902
\endref


\ref \key {Wa} \by L Washington \book Introduction to cyclotomic fields
\publ Springer--Verlag \yr 1982\endref

\ref \key {W} \by E Witten\paper Quantum field theory and the Jones
polynomial\jour Comm. Math. Phys.\vol 121 \yr 1989\pages 351--399
\endref

\endRefs
\bye